\documentclass[nospthms,envcountsect]{svmult} 
\usepackage{latexsym}
\usepackage{amsfonts}
\usepackage{amsmath}
\usepackage{amssymb}
\usepackage{amscd}
\usepackage{enumerate}
    \newtheorem{theorem}                    {Theorem}       [section]
    \newtheorem{lemma}      [theorem]       {Lemma}
    \newtheorem{corollary}  [theorem]       {Corollary}
    \newtheorem{proposition}[theorem]       {Proposition}
    
    \newtheorem{conjecture} [theorem]       {Conjecture}

\begin{document}
\catcode`@=11
\atdef@ I#1I#2I{\CD@check{I..I..I}{\llap{$\m@th\vcenter{\hbox
  {$\scriptstyle#1$}}$}
  \rlap{$\m@th\vcenter{\hbox{$\scriptstyle#2$}}$}&&}}
\atdef@ E#1E#2E{\ampersand@
  \ifCD@ \global\bigaw@\minCDarrowwidth \else \global\bigaw@\minaw@ \fi
  \setboxz@h{$\m@th\scriptstyle\;\;{#1}\;$}%
  \ifdim\wdz@>\bigaw@ \global\bigaw@\wdz@ \fi
  \@ifnotempty{#2}{\setbox@ne\hbox{$\m@th\scriptstyle\;\;{#2}\;$}%
    \ifdim\wd@ne>\bigaw@ \global\bigaw@\wd@ne \fi}%
  \ifCD@\enskip\fi
    \mathrel{\mathop{\hbox to\bigaw@{}}%
      \limits^{#1}\@ifnotempty{#2}{_{#2}}}%
  \ifCD@\enskip\fi \ampersand@}
\catcode`@=\active

\renewcommand{\labelenumi}{\alph{enumi})}
\newcommand{\chr}{\operatorname{char}}
\newcommand{\etale}{etale}
\newcommand{\isom}{\stackrel{\sim}{\longrightarrow}}
\newcommand{\Aut}{\operatorname{Aut}}
\newcommand{\Hom}{\operatorname{Hom}}
\newcommand{\End}{\operatorname{End}}
\newcommand{\HOM}{\operatorname{{\mathcal H{\frak{om}}}}}
\newcommand{\EXT}{\operatorname{\mathcal E{\frak xt}}}
\newcommand{\Tot}{\operatorname{Tot}}
\newcommand{\Cor}{\operatorname{Cor}}
\newcommand{\Ext}{\operatorname{Ext}}
\newcommand{\Tor}{\operatorname{Tor}}
\newcommand{\Gal}{\operatorname{Gal}}
\newcommand{\Pic}{\operatorname{Pic}}
\newcommand{\Spec}{\operatorname{Spec}}
\newcommand{\trdeg}{\operatorname{trdeg}}
\newcommand{\im}{\operatorname{im}}
\newcommand{\coim}{\operatorname{coim}}
\newcommand{\coker}{\operatorname{coker}}
\newcommand{\gr}{\operatorname{gr}}
\newcommand{\id}{\operatorname{id}}
\newcommand{\Br}{\operatorname{Br}}
\newcommand{\cd}{\operatorname{cd}}
\newcommand{\CH}{CH}
\newcommand{\Alb}{\operatorname{Alb}}
\renewcommand{\lim}{\operatornamewithlimits{lim}}
\newcommand{\colim}{\operatornamewithlimits{colim}}
\newcommand{\rk}{\operatorname{rank}}
\newcommand{\codim}{\operatorname{codim}}
\newcommand{\NS}{\operatorname{NS}}
\newcommand{\cone}{{\rm cone}}
\newcommand{\rank}{\operatorname{rank}}
\newcommand{\ord}{{\rm ord}}
\newcommand{\g}{{\cal G}}
\newcommand{\f}{{\cal F}}
\newcommand{\du}{{\cal D}}
\newcommand{\G}{{\mathbb G}}
\newcommand{\N}{{\mathbb N}}
\newcommand{\A}{{\mathbb A}}
\newcommand{\Z}{{{\mathbb Z}}}
\newcommand{\Q}{{{\mathbb Q}}}
\newcommand{\R}{{{\mathbb R}}}
\newcommand{\B}{{\mathbb Z}^c}
\renewcommand{\H}{{{\mathbb H}}}
\renewcommand{\P}{{{\mathbb P}}}
\newcommand{\F}{{{\mathbb F}}}
\newcommand{\m}{{\mathfrak m}}
\newcommand{\Sch}{{\text{\rm Sch}}}
\newcommand{\et}{{\text{\rm et}}}
\newcommand{\Zar}{{\text{\rm Zar}}}
\newcommand{\Nis}{{\text{\rm Nis}}}
\newcommand{\tr}{\operatorname{tr}}
\newcommand{\tor}{{\text{\rm tor}}}
\newcommand{\red}{{\text{\rm red}}}
\newcommand{\Div}{\operatorname{Div}}
\newcommand{\Ab}{{\text{\rm Ab}}}
\newcommand{\DD}{{\mathbb Z}^c}
\renewcommand{\div}{\operatorname{div}}
\newcommand{\corank}{\operatorname{corank}}
\renewcommand{\O}{{\cal O}}
\newcommand{\C}{{C}}
\newcommand{\p}{{\mathfrak p}}
\newcommand{\proof}{\noindent{\it Proof. }}
\newcommand{\proofend}{\hfill $\Box$ \\}
\newcommand{\rem}{\noindent {\it Remark. }}
\newcommand{\example}{\noindent {\bf Example. }}
\newcommand{\ar}{{\text{\rm ar}}}
\newcommand{\del}{{\delta}}
\title*{On Suslin's singular homology and cohomology}
\author{Thomas Geisser\thanks{Supported in part by NSF grant No.0901021}}
\institute{University of Southern California}

\maketitle

\vskip-1.5cm

\begin{abstract}
We study properties of Suslin homology and cohomology over non-algebraically
closed base fields, and their $p$-part in characteristic $p$.
In the second half we focus on finite fields, and consider finite
generation questions and connections to tamely ramified class
field theory.
\end{abstract}

\section{Introduction}
Suslin and Voevodsky defined Suslin homology (also called singular 
homology) $H_i^S(X,A)$ of a scheme of finite type over a field
$k$ with coefficients in an abelian group $A$ as 
$\Tor_i(\Cor_k(\Delta^*,X), A)$. Here $\Cor_k(\Delta^i,X)$ is 
the free abelian group generated by integral subschemes $Z$
of $\Delta^i\times X$ which are finite and surjective over $\Delta^i$,
and the differentials are given by alternating sums of pull-back maps 
along face maps.
Suslin cohomology $H^i_S(X,A)$ is defined to be 
$\Ext^i_\Ab(\Cor_k(\Delta^*,X),A)$. Suslin and Voevodsky showed
in \cite{susvoe98} that over a separably closed field in which $m$ 
is invertible, one has 
\begin{equation}\label{susvoemain}
H^i_S(X,\Z/m)\cong H^i_\et(X,\Z/m)
\end{equation}
(see \cite{ichdejong} for the case of fields of characteristic $p$).

In the first half of this paper we study both the situation 
that $m$ is a power of the characteristic of $k$, and that $k$
is not algebraically closed.
In the second half we focus on finite base fields and discuss a modified
version of Suslin homology, which is closely related to etale 
cohomology on the one hand,  but is also expected to be finitely
generated. Moreover, its zeroth homology is $\Z^{\pi_o(X)}$ and 
its first homology is expected to be an integral model of the abelianized 
tame fundamental group.

We start by discussing the $p$-part of Suslin homology over an 
algebraically closed field of characteristic $p$. 
We show that assuming resolution of 
singularities, the groups $H_i^S(X,\Z/p^r)$ are finite
abelian groups, and vanish outside the range $0\leq i \leq \dim X$.
Thus Suslin cohomology with finite coefficients
is etale cohomology away from the characteristic, but better behaved
than etale cohomology at the characteristic (for example, 
$H^1_\et({\mathbb A}^1,\Z/p)$ is not finite).
Moreover, Suslin homology is a birational invariant in the following 
strong sense:
If $X$ has a resolution of singularities $p:X'\to X$ which is an 
isomorphism outside of the open subset $U$, then 
$H_i^S(U,\Z/p^r)\cong H_i^S(X,\Z/p^r)$. 

Next we examine the situation over non-algebraically closed fields.
We redefine Suslin homology and cohomology by imposing 
Galois descent. Concretely, if $G_k$ is the absolute Galois
group of $k$, then we define Galois-Suslin homology to be
\begin{equation*}\label{galsushom}
H_i^{GS}(X,A) = H^{-i}R\Gamma(G_k,\Cor_{\bar k}(\Delta_{\bar k}^*,\bar X)),
\end{equation*}
and Galois-Suslin cohomology to be
\begin{equation*}\label{galsuscoh}
H^i_{GS}(X,A) = \Ext^i_{G_k}(\Cor_{\bar k}(\Delta_{\bar k}^*,\bar X), A).
\end{equation*}
Ideally one would like to define Galois-Suslin homology via Galois homology,
but we are not aware of such a theory. With rational coefficients,
the newly defined groups agree with the original groups. On the other hand,
with finite coefficients prime to the characteristic, the
proof of \eqref{susvoemain} in \cite{susvoe98} carries over to show that
$H^i_{GS}(X,\Z/m)\cong H^i_\et(X,\Z/m)$. As a corollary we obtain an 
isomorphism between $H_0^{GS}(X,\Z/m)$ and the abelianized fundamental group 
$\pi_1^{ab}(X)/m$ for any separated $X$ of finite type over a finite field. 

The second half of the paper focuses on the case of a finite base field.
We work under the assumption of resolution of singularities in order
to see the picture of the properties which can expected. 
The critical reader 
can view our statements as theorems for schemes of dimension at most three,
and conjectures in general.
A theorem of Jannsen-Saito \cite{uweshuji} can 
be generalized to show that Suslin homology and cohomology with finite
coefficients for any $X$ over a finite field is finite. 
Rationally, $H_0^S(X,\Q)\cong H^0_S(X,\Q)\cong \Q^{\pi_0(X)}$. Most other 
properties are equivalent to the following Conjecture $P_0$ considered in 
\cite{ichkato}: 
For $X$ smooth and proper over a finite field, $CH_0(X,i)$ is 
torsion for $i\not=0$. This is a particular case of Parshin's conjecture
that $K_i(X)$ is torsion for $i\not=0$. 
For example, Conjecture $P_0$ is equivalent to the vanishing of 
$H_i^S(X,\Q)$ for $i\not=0$ and all smooth $X$. For arbitrary $X$ 
of dimension $d$, Conjecture $P_0$ implies the
vanishing of $H_i^S(X,\Q)$ outside of the range $0\leq i\leq d$ 
and its finite dimensionality in this range. 
Combining the torsion
and rational case, we show that $H_i^S(X,\Z)$ and $H^i_S(X,\Z)$ 
are finitely generated for all $X$ if and only if Conjecture $P_0$ holds.

Over a finite field and with integral coefficients, it is more
natural to impose descent by the Weil group $G$ generated by
the Frobenius endomorphism instead of the Galois group
\cite{lichtenbaumweilI, ichweilI, ichduke, ichkato}.
We define arithmetic homology 
$$H_i^\ar(X,A)=\Tor_i^G(\Cor_{\bar k}(\Delta_{\bar k}^*,\bar X),A)$$ 
and arithmetic cohomology 
$$H^i_\ar(X,\Z)=\Ext^i_G(\Cor_{\bar k}(\Delta_{\bar k}^*,\bar X),\Z).$$ 
We chose the notation to distinguish the groups from weight homology
and Weil-etale cohomology considered elsewhere.  
We show that $H_0^\ar(X,\Z)\cong H^0_\ar(X,\Z)\cong \Z^{\pi_0(X)}$
and that arithmetic homology and cohomology lie in long exact
sequences with Galois-Suslin homology and cohomology, respectively. 
They are finitely generated abelian groups if and only if
Conjecture $P_0$ holds. 

The difference between arithmetic and Suslin homology
is measured by a third theory, which we call Kato-Suslin
homology, and which is defined as 
$H_i^{KS}(X,A)=H_i(\Cor_{\bar k}(\Delta^*_{\bar k},\bar X)_G)$.
By definition there is a long exact sequence
$$ \cdots \to H_i^S(X,A) \to  H_{i+1}^\ar(X,A) \to 
H_{i+1}^{KS}(X,A) \to H_{i-1}^{S}(X,A) \to \cdots .$$
It follows that $H_0^{KS}(X,A)=\Z^{\pi_0(X)}$ for any $X$. 
As a generalization of the integral version \cite{ichkato}
of Kato's conjecture \cite{kato}, we propose 

\begin{conjecture} The groups $H_i^{KS}(X,A)$ vanish for all $i>0$
and all smooth $X$.
\end{conjecture} 
 
Equivalently that there are short exact sequences 
$$0\to H_{i+1}^S(\bar X,\Z)_G \to H_i^S(X,\Z) \to H_i^S(\bar X,\Z)^G\to 0$$ 
for all $i\geq 0$ and all smooth $X$.
We show that this conjecture, too, is equivalent to Conjecture $P_0$.
This leads us to conjecture on class field theory:

\begin{conjecture} For every $X$ separated and of finite type over $\F_q$, 
there is a canonical injection
$$H_1^\ar(X,\Z)\to \pi_1^t(X)^{ab}$$ 
with dense image.
\end{conjecture} 
 
It might even be true that the relative group 
$H_1^\ar(X,\Z)^\circ := \ker (H_1^\ar(X,\Z)\to \Z^{\pi_0(X)})$
is isomorphic to the geometric part of the abelianized fundamental group
defined in SGA 3X\S 6. 
To support our conjecture, we note that the generalized Kato conjecture 
above implies $H_0^S(X,\Z)\cong H_1^\ar(X,\Z)$ for smooth $X$, so that 
in this case our conjecture becomes a theorem of Schmidt-Spiess \cite{ss}.
In addition, we show (independent of any conjectures) 

\begin{proposition}
If $\frac{1}{l}\in \F_q$, then 
$H_1^\ar(X,\Z)^{\wedge l}\cong \pi_1^t(X)^{ab}(l)$ for arbitrary $X$.
\end{proposition}

In particular, the conjectured finite generation of $H_1^\ar(X,\Z)$ implies 
the conjecture away from the characteristic. 
We also give a conditional result at the characteristic.

Notation: In this paper, scheme over a field $k$ means separated scheme
of finite type over $k$. 

We thank Uwe Jannsen for interesting discussions
related to the subject of this paper, and Shuji Saito and Takeshi Saito
for helpful comments during a series of lectures I gave on the topic
of this paper at Tokyo University.

\section{Motivic homology}
Suslin homology $H_i^S(X,\Z)$ of a scheme $X$ over a
field $k$ is defined as the homology of the global sections
$C_*^X(k)$ of the complex of etale sheaves
$C_*^X(-)=\Cor_k(-\times \Delta^*,X)$. Here 
$\Cor_k(U,X)$ is the group 
of universal relative cycles of $U\times Y/U$ \cite{susvoerel}.
If $U$ is smooth, then $\Cor_k(U,X)$
is the free abelian group generated by closed irreducible subschemes
of $U\times X$ which are finite and surjective over a connected
component of $U$.

More generally \cite{friedvoe},
motivic homology of weight $n$ are the extension groups in Voevodsky's category of geometrical mixed motives
$$ H_i(X,\Z(n))=\Hom_{DM^-_k}(\Z(n)[i],M(X)),$$
and are isomorphic to
$$H_i(X,\Z(n))=\begin{cases}
H^{2n-i}_{(0)}({\mathbb A}^n,C_*^X)&n\geq 0\\
H_{i-2n-1}(C_*\big(\frac{c_0(X\times ({\mathbb A}^n-\{0\}))}
{c_0(X\times \{1\})}\big)(k)) &n< 0.
\end{cases}$$
Here cohomology is taken for the Nisnevich topology.
There is an obvious version with coefficients.
Motivic homology is a covariant functor on the category of schemes
of finite type over $k$, and has the following
additional properties, see \cite{friedvoe} (the final three properties
require resolution of singularities)

\begin{enumerate}
\item It is homotopy invariant.
\item It satisfies a projective bundle formula
$$ H_i(X\times {\mathbb P}^1,\Z(n))=
H_i(X,\Z(n))\oplus H_{i-2}(X,\Z(n-1)).$$
\item There is a Mayer-Vietoris long exact sequence for open covers.
\item  Given an abstract blow-up square
$$\begin{CD}
Z'@>>> X'\\
@VVV @VVV \\
Z@>>> X
\end{CD}$$
there is a long exact sequence
\begin{multline}\label{blowups}
\cdots\to H_{i+1}(X,\Z(n))\to
H_i(Z',\Z(n))\to \\
H_i(X',\Z(n))\oplus H_i(Z,\Z(n))\to H_i(X,\Z(n))\to \cdots
\end{multline}
\item If $X$ is proper, then motivic homology agrees with higher
Chow groups, $H_i(X,\Z(n))\cong CH_n(X,i-2n)$.
\item If $X$ is smooth of pure dimension $d$, then motivic homology
agrees with motivic cohomology with compact support,
$$H_i(X,\Z(n))\cong H^{2d-i}_c(X,\Z(d-n)).$$
In particular, if $Z$ is
a closed subscheme of a smooth scheme $X$ of pure dimension $d$,
then we have a long exact sequence
\begin{equation}\label{localizations}
\cdots \to H_i(U,\Z(n)) \to H_i(X,\Z(n)) \to
H^{2d-i}_c(Z,\Z(d-n))\to\cdots .
\end{equation}
\end{enumerate}

In order to remove the hypothesis on resolution of singularities,
it would be sufficient to find a proof of Theorem 5.5(2) of 
Friedlander-Voevodsky \cite{friedvoe} that does
not require resolution of singularities.
For all arguments in this paper (except the $p$-part of the Kato
conjecture) the sequences \eqref{blowups} and \eqref{localizations}
and the existence of a smooth and 
proper model for every function field are sufficient.

\subsection{Suslin cohomology}
Suslin cohomology is by definition the dual of Suslin homology,
i.e. for an abelian group $A$ it is defined as
$$ H^i_S(X,A):= \Ext^i_\Ab(\C_*^X(k),A).$$
We have $H^i_S(X,\Q/\Z)\cong \Hom(H_i^S(X,\Z),\Q/\Z)$, and
a short exact sequence of abelian groups gives a long exact sequence of
cohomology groups, in particular long exact sequences
\begin{equation}\label{coeffseq}
\cdots \to H^i_S(X,\Z)\to
H^i_S(X,\Z)\to H^i_S(X,\Z/m)\to H^{i+1}_S(X,\Z)\to \cdots.
\end{equation}
and
$$\cdots \to H^{i-1}_S(X,\Q/\Z)\to H^i_S(X,\Z)\to
H^i_S(X,\Q)\to H^i_S(X,\Q/\Z)\to \cdots.$$
Consequently, $H^i_S(X,\Z)_\Q\cong H^i_S(X,\Q)$ if
Suslin-homology is finitely generated.
If $A$ is a ring, then $H^i_S(X,A)\cong \Ext_A^i(\C_*^X(k)\otimes A,A)$,
and we get a spectral sequence
\begin{equation}\label{extss}
E_2^{s,t}= \Ext^s_A(H_t^S(X,A),A)\Rightarrow H^{s+t}_S(X,A).
\end{equation}
In particular, there are perfect pairings
\begin{align*}
H_i^S(X,\Q)\times H^i_S(X,\Q) &\to \Q\\
H_i^S(X,\Z/m)\times H^i_S(X,\Z/m) &\to \Z/m.
\end{align*}

\begin{lemma}\label{duall1}
There are natural pairings
$$ H^i_S(X,\Z)/\tor \times H_i^S(X,\Z)/\tor \to \Z$$
and
$$H^i_S(X,\Z)_\tor \times H_{i-1}^S(X,\Z)_\tor\to \Q/\Z.$$
\end{lemma}

\proof
The spectral sequence \eqref{extss} gives a short exact sequence
\begin{equation}\label{cefff}
0\to \Ext^1(H_{i-1}^S(X,\Z),\Z) \to H^i_S(X,\Z)\to
\Hom(H_i^S(X,\Z),\Z)\to 0.
\end{equation}
The resulting map $H^i_S(X,\Z)/\tor \twoheadrightarrow \Hom(H_i^S(X,\Z),\Z)$
induces the first pairing. Since $\Hom(H_i^S(X,\Z),\Z)$
is torsion free, we obtain the map 
\begin{multline*}
H^i_S(X,\Z)_\tor \hookrightarrow 
\Ext^1(H_{i-1}^S(X,\Z),\Z)\twoheadrightarrow \\
\Ext^1(H_{i-1}^S(X,\Z)_\tor ,\Z) \stackrel{\sim}{\leftarrow} 
\Hom(H_{i-1}^S(X,\Z)_\tor,\Q/\Z)
\end{multline*}
for the second pairing.
\proofend

\subsection{Comparison to motivic cohomology}
Recall that in the category $DM^-_k$ of bounded above complexes
of homotopy invariant Nisnevich sheaves with transfers,
the motive $M(X)$ of $X$ is the complex of presheaves with
transfers $\C_*^X$.
Since a field has no higher Nisnevich cohomology,
taking global sections over $k$ induces a canonical map 
$$\Hom_{DM^-_k}(M(X),A) \to R\Hom_{\Ab}(C_*^X(k),R\Gamma(k,A)),$$
hence a natural map
\begin{equation}\label{canmap}
H^i_M(X,A)\to H^i_S(X,A).
\end{equation}
Even though the cohomology groups do
not depend on the base field, the map does. For example, if $L/k$
is an extension of degree $d$, then the diagram of groups isomorphic 
to $\Z$,
$$\begin{CD}
H^0_M(\Spec k,\Z)@= H^0_S(\Spec k,\Z)\\
@| @VV\times d V\\
H^0_M(\Spec L,\Z)@>>> H^0_S(\Spec L,\Z)
\end{CD}$$
shows that the lower horizontal map is
multiplication by $d$.  We will see
below that conjecturally \eqref{canmap} is a map between finitely
generated groups which is rationally an isomorphism, and one
might ask if its Euler characteristic has any interpretation.

\section{The mod $p$ Suslin homology in characteristic $p$}
We examine the $p$-part of Suslin homology
in characteristic $p$. We assume
that $k$ is perfect and resolution of singularities exists
over $k$ in order to obtain stronger results.
We first give an auxiliary result on motivic cohomology
with compact support:

\begin{proposition}\label{theelemma}
Let $d=\dim X$.

a) We have $H^i_c(X,\Z/p^r(n))=0$ for $n>d$.

b) If $k$ is algebraically closed, then $H^i_c(X,\Z/p^r(d))$ is finite,
$H^i_c(X,\Q_p/\Z_p(d))$ is of cofinite type, and the groups vanish
unless $d\leq i\leq 2d$.
\end{proposition}

\proof
By induction on the dimension and the localization sequence,
the statement for $X$ and a dense open subset of $X$
are equivalent. Hence replacing $X$ by a smooth subscheme and then
by a smooth and proper model, we can assume that $X$ is smooth
and proper. Then a) follows from \cite{ichmarcI}.  If $k$ is
algebraically closed, then 
$$ H^i(X,\Z/p(d))\cong H^{i-d}(X_{Nis},\nu^d)\cong
H^{i-d}(X_\et,\nu^d),$$
by \cite{ichmarcI} and \cite{katokuzumaki}. That the latter
group is finite and of cofinite type, respectively, can be derived
from \cite[Thm.1.11]{milnevalues}, and it vanishes outside of the given
range by reasons of cohomological dimension.
\proofend

\begin{theorem}\label{ppart}
Let $X$ be separated and of finite type over $k$.

a) The groups $H_i(X,\Z/p^r(n))$ vanish for all $n<0$.

b) If $k$ is algebraically closed, then the groups 
$H_i^S(X,\Z/p^r)$ are finite,
the groups $H_i^S(X,\Q_p/\Z_p)$ are of cofinite type, and both vanish
unless $0\leq i\leq d$.
\end{theorem}

\proof
If $X$ is smooth, then $H_i(X,\Z/p^r(n))\cong H^{2d-i}_c(X,\Z/p^r(d-n))$
and we conclude by the Proposition.
In general, we can assume by \eqref{blowups} and induction on the
number of irreducible components that $U$ is integral. Proceeding
by induction on the dimension, we choose a resolution of singularities
$X'$ of $X$, let $Z$ be the closed subscheme of $X$ where the map 
$X'\to X$ is not an isomorphism, and let $Z'=Z\times_XX'$.
Then we conclude by the sequence \eqref{blowups}.
\proofend

\example If $X'$ is the blow up of a smooth scheme $X$ in a smooth 
subscheme $Z$, then the strict transform $Z'=X'\times_XZ$ is a projective
bundle over $Z$, hence by the projective bundle formula
$H_i^S(Z,\Z/p^r)\cong H_i^S(Z',\Z/p^r)$ and 
$H_i^S(X,\Z/p^r)\cong H_i^S(X',\Z/p^r)$. More generally, we have

\begin{proposition}\label{biratt}
Assume $X$ has a desingularization $p:X'\to X$ which is 
an isomorphism outside of the dense open subset $U$. Then 
$H_i^S(U,\Z/p^r)\cong  H_i^S(X,\Z/p^r)$.
In particular, the $p$-part of Suslin homology is a birational invariant.
\end{proposition}

The hypothesis is satisfied if $X$ is smooth, or if $U$ contains all 
singular points of $X$ and
a resolution of singularities exists which is an isomorphism outside
of the singular points.

\medskip

\proof
If $X$ is smooth, then this follows from Proposition \ref{theelemma} 
and the localization sequence \eqref{localizations}.
In general, let $Z$ be the set of points where $p$ is not an isomorphism,
and consider the cartesian diagram
$$\begin{CD}
Z' @>>> U' @>>> X'\\
@VVV @VVV @VVV \\
Z@>>> U @>>> X.
\end{CD}$$
Comparing long exact sequence \eqref{blowups} of the left and outer
squares, 
$$\begin{CD}
\to H_i^S(Z',\Z/p^r)@>>>  H_i^S(U',\Z/p^r)\oplus H_i^S(Z,\Z/p^r)@>>>
H_i^S(U,\Z/p^r)\to \\
@| @|@VVV \\
\to H_i^S(Z',\Z/p^r)@>>>  H_i^S(X',\Z/p^r)\oplus H_i^S(Z,\Z/p^r)@>>>
H_i^S(X,\Z/p^r)\to
\end{CD}$$
we see that
$H_i^S(U',\Z/p^r)\cong  H_i^S(X',\Z/p^r)$ implies
$H_i^S(U,\Z/p^r)\cong  H_i^S(X,\Z/p^r)$.
\proofend

\example If $X$ is a node, then the blow-up sequence gives
$H_i^S(X,\Z/p^r)\cong H_{i-1}^S(k,\Z/p^r)\oplus H_i^S(k,\Z/p^r)$,
which is $\Z/p^r$ for $i=0,1$ and vanishes otherwise. 
Reid constructed a normal surface with a singular point
whose blow-up is a node, showing that
the $p$-part of Suslin homology is not a birational invariant for normal
schemes.

\begin{corollary}
The higher Chow groups $CH_0(X,i,\Z/p^r)$ and the logarithmic de Rham-Witt
cohomology groups $H^i(X_\et,\nu_r^d)$, for $d=\dim X$,  
are birational invariants.
\end{corollary}

\proof
Suslin homology agrees with higher Chow groups for proper $X$,
and with motivic cohomology for smooth and proper $X$.
\proofend

Note that integrally $CH_0(X)$ is a birational invariant, but
the higher Chow groups $CH_0(X,i)$ are generally not.

Suslin and Voevodsky \cite[Thm.3.1]{susvoe98} show that for a smooth
compactification $\bar X$ of the smooth curve $X$, 
$H_0^S(X,\Z)$ is isomorphic to the relative Picard group $\Pic(\bar X,Y)$
and that all higher Suslin homology groups vanish.
Proposition \ref{biratt} implies that the
kernel and cokernel of $\Pic(\bar X,Y)\to \Pic(\bar X)$
are uniquely $p$-divisible.
Given $U$ with compactification $j:U\to X$, the
normalization $X^\sim$ of $X$ in $U$ is the affine bundle defined
by the integral closure of $\O_X$ in $j_*\O_U$. We call $X$ normal
in $U$ if $X^\sim\to X$ is an isomorphism.

\begin{proposition} If $X$ is normal in the curve $U$, then 
$H_i^S(U,\Z/p)\cong H_i^S(X,\Z/p)$.
\end{proposition}

\proof
This follows by applying the argument of Proposition \ref{biratt}
to $X'$ the normalization of $X$, $Z$ the closed subset
where $X'\to X$ is not an isomorphism, $Z'=X'\times_XZ$
and $U'=X'\times_XU$. Since $X$ is normal in $U$, we have
$Z\subseteq U$ and $Z'\subseteq U'$.
\proofend

\section{Galois properties}
Suslin homology is covariant, i.e. a separated map $f:X\to Y$ of 
finite type induces
a map $f_*:\Cor_k(T,X)\to \Cor_k(T,Y)$ by sending a closed
irreducible subscheme $Z$ of $T\times X$, finite over $T$,
to the subscheme $[k(Z):k(f(Z))]\cdot f(Z)$ (note that $f(Z)$
is closed in $T\times Y$ and finite over $T$). On the other
hand, Suslin homology
is contravariant for finite flat maps $f:X\to Y$, because $f$
induces a map $f^*:\Cor_k(T,Y)\to \Cor_k(T,X)$ by composition with the
graph of $f$ in $\Cor_k(Y,X)$ (note that the graph is a universal
relative cycle in the sense of \cite{susvoerel}).
We examine the properties of Suslin homology under change of base-fields.

\begin{lemma}
Let $L/k$ be a finite extension of fields, $X$ a scheme over $k$
and $Y$ a scheme over $L$. Then 
$ \Cor_L(Y,X_L)= \Cor_k(Y,X)$ and if $X$ is smooth, then
$ \Cor_L(X_L,Y)= \Cor_k(X,Y)$. 
In particular, Suslin homology does not depend on the base field.
\end{lemma}

\proof
The first statement follows because $Y\times_LX_L\cong Y\times_kX$. The
second statement follows because
the map $X_L\to X$ is finite and separated, hence a closed subscheme of
$X_L\times_LY \cong X\times_kY$ is finite and surjective over
$X_L$ if and only if it is finite and surjective over $X$.
\proofend

Given a scheme over $k$, the graph of the projection $X_L\to X$
in $X_L\times X$ gives elements $\Gamma_X\in \Cor_k(X,X_L)$
and $\Gamma_X^t\in \Cor_k(X_L,X)$.

\subsection{Covariance}
\begin{lemma}
a) If $X$ and $Y$ are separated schemes of finite type over $k$, then 
the two maps
$$\Cor_L(X_L,Y_L)\to \Cor_k(X,Y)$$
induced by composition and precomposition, respectively,
with $\Gamma_Y^t$ and $\Gamma_X$ agree. Both maps send a generator
$Z\subseteq X_L\times_k Y\cong X\times_k Y_L$ to its image in $X\times Y$
with multiplicity $[k(Z):k(f(Z))]$, a divisor of $[L:k]$.

b) If $F/k$ is an infinite algebraic extension, then 
$\lim_{L/k} \Cor_L(X_L,Y_L)=0$.
\end{lemma}

\proof
The first part is easy. If $Z$ is of finite type over $k$, 
then $k(Z)$ is a finitely generated field extension of $k$. 
For every component $Z_i$ of $Z_F$, we 
obtain a map $F\to F\otimes_k k(Z) \to k(Z_i)$, and
since $F$ is not finitely generated over $k$, neither is $k(Z_i)$.
Hence going up the tower of finite
extensions $L/k$ in $F$, the degree of $[k(W_L):k(Z)]$, for $W_L$ the
component of $Z_L$ corresponding to $Z_i$, goes to infinity. 
\proofend


\subsection{Contravariance}

\begin{lemma}\label{contra}
a) If $X$ and $Y$ are schemes over $k$, then the two maps
$$\Cor_k(X,Y)\to \Cor_L(X_L,Y_L)$$
induced by composition and precomposition, respectively with $\Gamma_Y$ and
$\Gamma_X^t$ agree. Both maps send a generator
$Z\subseteq X\times Y$ to the cycle associated to
$Z_L\subseteq X\times_k Y_L\cong X_L\times_k Y$. 
If $L/k$ is separable, this is a sum of the integral
subschemes lying over $Z$ with multiplicity one.
If $L/k$ is Galois with group $G$, then the maps induce an isomorphism
$$\Cor_k(X,Y)\cong  \Cor_L(X_L,Y_L)^G.$$

b) Varying $L$, $\Cor_L(X_L,Y_L)$ forms a \etale sheaf on
$\Spec k$ with stalk
$M=\colim_L \Cor_L(X_L,Y_L)\cong \Cor_{\bar k}(X_{\bar k},Y_{\bar k})$,
where $L$ runs through the finite extensions of $k$ in an
algebraic closure $\bar k$ of $k$. In particular,
$\Cor_L(X_L,Y_L)\cong M^{\Gal(\bar k/L)}$.
\end{lemma}

\proof
Again, the first part is easy. If $L/k$ is separable, 
$Z_L$ is finite and \etale over $Z$, 
hence $Z_L\cong \sum_i Z_i$, a finite sum of the integral
cycles lying over $Z$ with multiplicity one each.
If $L/k$ is moreover Galois, then $\Cor_k(X,Y)\cong \Cor_L(X_L,Y_L)^G$ and
$\Cor_{\bar k}(X_{\bar k},Y_{\bar k})\cong \colim \Cor_L(X_L,Y_L)$
by EGA IV Thm. 8.10.5.
\proofend

The proposition suggests to work with the complex $C_*^X $ of 
etale sheaves on $\Spec k$ given by 
$$\C_*^X(L):=\Cor_L(\Delta^*_L,X_L)\cong \Cor_k(\Delta^*_L,X).$$

\begin{corollary}\label{descentss}
We have $H_i^S(\bar X,A)\cong \colim_L H_i^S(X_L,A)$, and
there is a spectral sequence
$$E_2^{s,t}=  {\lim}^s H^t_S(X_L,A) \Rightarrow
H^{s+t}_S(\bar X,A).$$
The maps in the direct and inverse system are induced by contravariant
functoriality of Suslin homology for finite flat maps.
\end{corollary}

\proof
This follows from the quasi-isomorphisms
$$ R\Hom_\Ab(\C_*^X(\bar k),\Z)
\cong R\Hom_\Ab(\colim_L \C_*^X(L),\Z)
\cong R\lim_L R\Hom_\Ab(\C_*^X(L),\Z).$$
\proofend

\subsection{Coinvariants}
If $G_k$ is the absolute Galois group of $k$, then 
$\Cor_{\bar k}(\bar X,\bar Y)_{G_k}$ can be identified
with $\Cor_{k}(X,Y)$ by associating orbits of points
of $\bar X\times_{\bar k}\bar Y$ with their image in $X\times_k Y$.
However, this identification is neither compatible with covariant nor with
contravariant functoriality, and in particular not with the differentials
in the complex $C_*^X(k)$.
But the obstruction is torsion, and we can remedy this problem by
tensoring with $\Q$: Define an isomorphism
$$\tau:(\Cor_{\bar k}(\bar X,\bar Y)_\Q)_{G_k}
\to \Cor_{k}(X,Y)_\Q.$$
as follows. A generator $1_{\bar Z}$ corresponding to the closed
irreducible subscheme $\bar Z\subseteq \bar X\times\bar Y$ is sent to
$\frac{1}{g_Z} 1_Z$, where $Z$ is the image of $\bar Z$ in $X\times Y$
and $g$ the number of irreducible components of $Z\times_k\bar k$,
i.e. $g_Z$ is the size of the Galois orbit of $\bar Z$.

\begin{lemma}\label{ratisom}
The isomorphism $\tau$ is functorial in both variables,
hence it induces an isomorphism of complexes
$$( C_*^X(\bar k)_\Q )_{G_k} \cong C_*^X(k)_\Q.$$
\end{lemma}

\proof
This can be proved by direct verification. 
We give an alternate proof. Consider the composition
$$\Cor_{k}(X,Y)\to \Cor_{\bar k}(\bar X,\bar Y)^{G_k}\to
\Cor_{\bar k}(\bar X,\bar Y)_{G_k}\xrightarrow{\tau} \Cor_{k}(X,Y)_\Q.$$
The middle map is induced by the identity, and is
multiplication by $g_Z$ on the component corresponding to $Z$.
All maps are isomorphisms upon tensoring with $\Q$. The first map, the  
second map, and the composition are functorial, hence so is the $\tau$.
\proofend

\section{Etale theory}
Let $\bar k$ be the algebraic closure of $k$ with Galois group $G_k$,
and let $A$ be a continuous $G_k$-module.
Then $\C_*^X(\bar k)\otimes A$
is a complex of continuous $G_k$-modules, and if $k$ has 
finite cohomological dimension we define Galois-Suslin homology to be
$$H_i^{GS}(X,A)= H^{-i}R\Gamma(G_k, \C_*^X(\bar k)\otimes A).$$
By construction, there is a spectral sequence
$$E^2_{s,t}=H^{-s}(G_k,H_t^S(\bar X,A))\Rightarrow
H_{s+t}^{GS}(X,A).$$
The case $X=\Spec k$ shows that Suslin homology does not agree with
\etale Suslin homology, i.e. Suslin homology does not have Galois descent.
We define \etale Suslin cohomology to be
\begin{equation}\label{suscohdef}
H^i_{GS}(X,A)= \Ext^i_{G_k}(\C_*^X(\bar k),A).
\end{equation}
This agrees with the old definition if $k$ is algebraically closed.
Let $\tau_*$ be the functor from $G_k$-modules to continuous
$G_k$-modules which sends $M$ to $\colim_L M^{G_L}$, where 
$L$ runs through the finite extensions of $k$. 
It is easy to see that $R^i\tau_* M= \colim_H H^i(H,M)$.

\begin{lemma}\label{suscohss}
We have  $H^i_{GS}(X,A)=H^i R\Gamma_{G_k}R\tau_*\Hom_\Ab(\C_*^X(\bar k), A)$.
In particular, there is a spectral sequence
\begin{equation}
E_2^{s,t}= H^s(G_k,R^t\tau_* \Hom_\Ab(\C_*^X(\bar k), A))\Rightarrow
H^{s+t}_{GS}(X,A).
\end{equation}
\end{lemma}

\proof
This is \cite[Ex. 0.8]{adt}. Since $C_*^X(\bar k)$ is a
complex of free $\Z$-modules, $\Hom_\Ab(\C_*^X(\bar k), -)$ is
exact and preserves injectives. Hence the derived functor of
$\tau_*\Hom_\Ab(\C_*^X(\bar k), -)$ is $R^t\tau_* $ applied
to $\Hom_\Ab(\C_*^X(\bar k), -)$.
\proofend

\begin{proposition}\label{Qversion}
If $A$ is a $\Q$-vector space with trivial $G_k$-action, then
\begin{align*}
H_i^{GS}(X,A)&\cong H_i^S(X,A)\\
H^i_{GS}(X,A)&\cong H^i_S(X,A).
\end{align*}
\end{proposition}

\proof
Since $-\otimes A$ is exact, 
$H_i^S(X,A)= H_i(C_*^X(\bar k)^{G_k}\otimes A)$
as well as $H_i^{GS}(X,A)=H_i((C_*^X(\bar k)\otimes A )^{G_k})$ 
are isomorphic to the homology of the kernel of the map of complexes
$$C_*^X(\bar k)\otimes A \stackrel{\varphi-1}{\to}C_*^X(\bar k)\otimes A.$$
Since higher Galois cohomology is torsion, we have
$R^t\tau_*\Hom(\C_i^X(\bar k), A)=0$ for $t>0$, and
$H^s(G_k,\tau_* \Hom(\C_*^X(\bar k), A))=0$ for $s>0$. Hence
$H^i_{GS}(X,A)$ is isomorphic to the $i$th cohomology of
$$
\Hom_{G_k}(\C^X_*(\bar k),A)
\cong \Hom_\Ab(\C^X_*(\bar k)_{G_k},A)
\cong \Hom_\Ab(\C_*^X(k),A).
$$
The latter equality follows with Lemma \ref{ratisom}.
\proofend

\begin{theorem}\label{svduality} 
If $m$ is invertible in $k$ and $A$ is a finitely generated 
$m$-torsion $G_k$-module, then
$$ H^i_{GS}(X,A)\cong H^i_{\et}(X,A).$$
\end{theorem}

\proof
This follows with the argument of Suslin-Voevodsky \cite{susvoe98}.
Indeed, let $f:(Sch/k)_h \to Et_k$ be the canonical map
from the large site with the h-topology of $k$ to the small etale site of $k$. 
Clearly $f_*f^*\f\cong \f$, and the proof of Thm.4.5 in loc.cit. shows that 
the cokernel of the injection $f^*f_* \f \to \f$ is uniquely $m$-divisible,
for any homotopy invariant presheaf with transfers 
(like, for example, $C_i^X: U\mapsto \Cor_k(U\times \Delta^i,X)$). Hence 
$$\Ext^i_h(\f_h^\sim,f^*A)\cong \Ext^i_h(f^*f_*\f_h^\sim,f^*A)\cong
\Ext^i_{Et_k}(f_*\f^\sim_h,A)\cong \Ext^i_{G_k}(\f(\bar k),A).$$
Then the argument of section 7 in loc.cit. together with 
Theorem 6.7 can be descended from the algebraic closure of $k$ to $k$.
\proofend

\subsection{Duality results}
Duality results for the Galois cohomology of a field $k$
lead via theorem \ref{svduality} to duality results between
Galois-Suslin homology and cohomology over $k$.

\begin{theorem}\label{ffduality}
Let $k$ be a finite field, $A$ a finite $G_k$-module, and
$A^*=\Hom(A,\Q/\Z)$.
Then there is a perfect pairing of finite groups
$$ H_{i-1}^{GS}(X,A) \times H^i_{GS}(X,A^*) \to \Q/\Z.$$
\end{theorem}

\proof
According to \cite[Example 1.10]{adt} we have
$$\Ext^r_{G_k}(M,\Q/\Z)\cong \Ext^{r+1}_{G_k}(M,\Z)
\cong H^{1-r}(G_k,M)^*.$$
Hence
\begin{multline*}
\Ext^r_{G_k}(C_*^X(\bar k),\Hom(A,\Q/\Z))\cong
\Ext^r_{G_k}(C_*^X(\bar k)\otimes A,\Q/\Z)\cong \\
H^{1-r}(G_k,C_*^X(\bar k)\otimes A)^*
\cong H_{r-1}^{GS}(X,A)^*.
\end{multline*}
\proofend

The case of non-torsion sheaves is discussed below.

\begin{theorem}
Let $k$ be a local field with finite residue field and separable
closure $k^s$. For a finite $G_k$-module $A$ let 
$A^D=\Hom(A,(k^s)^\times)$. 
Then we have isomorphisms
$$H^i_{GS}(X,A^D) \cong\Hom(H_{i-2}^{GS}(X,A) ,\Q/\Z).$$
\end{theorem}

\proof
According to \cite[Thm.2.1]{adt} we have
$$\Ext^r_{G_k}(M,(k^s)^ \times)\cong H^{2-r}(G_k,M)^*$$
for every finite $G_k$-module $M$.
The rest of the proof is the same as above.
\proofend

\section{Finite base fields}
From now on we fix a finite field $\F_q$ with algebraic
closure $\bar \F_q$. To obtain the following results, we 
assume resolution of singularities. This is
needed to use the sequences \eqref{blowups} and \eqref{localizations}
to reduce to the smooth
and projective case on the one hand, and the proof of 
Jannsen-Saito \cite{uweshuji} of the Kato conjecture on the other hand
(however, Kerz and Saito announced a proof of the prime to $p$-part
of the Kato conjecture which does not require resolution of 
singularities).
The critical reader is invited to view the following results as 
conjectures which are theorems in dimension at most $3$.
 
We first present results on finite generation in the spirit
of \cite{uweshuji} and \cite{ichkato}.

\begin{theorem}\label{jsfingen}
For any $X/\F_q$ and any integer $m$, 
the groups $H_i^S(X,\Z/m)$ and $H^i_S(X,\Z/m)$ are finitely generated.
\end{theorem}

\proof
It suffices to consider the case of homology.
If $X$ is smooth and proper of dimension $d$, then 
$H_i^S(X,\Z/m)\cong CH_0(X,i,\Z/m)\cong H^{2d-i}_c(X,\Z/m(d))$,
and the result follows from work of Jannsen-Saito \cite{uweshuji}.
The usual devisage then shows that $H^j_c(X,\Z/m(d))$ is finite
for all $X$ and $d\geq \dim X$, hence $H_i^S(X,\Z/m)$ is finite
for smooth $X$. Finally, one proceeds by induction on the dimension
of $X$ with the blow-up long-exact sequence to reduce to the case
$X$ smooth.
\proofend

\subsection{Rational Suslin-homology}
We have the following unconditional result:

\begin{theorem}\label{ratdegree0}
For every connected $X$, the map $H_0^S(X,\Q)\to H_0^S(\F_q,\Q)\cong \Q$
is an isomorphism.
\end{theorem}

\proof
By induction on the number of irreducible components and \eqref{blowups}
we can first assume that $X$ is 
irreducible and then reduce to the situation where $X$ is smooth. 
In this case, we use \eqref{localizations} and the following 
Proposition to reduce to the smooth and proper case, where
$H_0^S(X,\Q)=CH_0(X)_\Q \cong CH_0(\F_q)_\Q$.
\proofend

\begin{proposition}
If $n>\dim X$, then $H^i_c(X,\Q(n))=0$
for $i\geq n+\dim X$. 
\end{proposition}

\proof
By induction on the dimension and the localization sequence
for motivic cohomology with compact support one sees that the
statement for $X$ and a dense open subscheme of $X$ are equivalent.
Hence we can assume that $X$ is smooth and proper of dimension $d$.
Comparing to higher Chow groups, one sees that this vanishes for 
$i>d+n$ for dimension (of cycles) reasons.
For $i=d+n$, we obtain from the niveau spectral sequence a surjection
$$ \bigoplus_{X_{(0)}} H^{n-d}_M(k(x),\Q(n-d))
\twoheadrightarrow H^{d+n}_M(X,\Q(n)).$$
But the summands vanish for $n>d$ because higher Milnor $K$-theory of
finite fields is torsion. 
\proofend

By definition, the groups $H_i(X,\Q(n))$ vanish for $i<n$.
We will consider the following conjecture $P_n$ of \cite{ichparshin}:

\medskip

\noindent{\bf Conjecture $P_n$:}
{\it For all smooth and projective schemes $X$ over the finite field $\F_q$,
the groups $H_i(X,\Q(n))$ vanish for $i\not=2n$.}

\medskip

This is a special case of Parshin's conjecture: If $X$ is smooth
and projective of dimension $d$, then
$$ H_i(X,\Q(n))\cong H^{2d-i}_M(X,\Q(d-n))\cong K_{i-2n}(X)^{(d-n)}$$
and, according to Parshin's conjecture, the latter $K$-group vanishes for 
$i\not=2n$. 
By the projective bundle formula, $P_n$ implies $P_{n-1}$.

\begin{proposition}\label{ratsituation}
a) Let $U$ be a curve. Then $H_i^S(U,\Q)\cong H_i^S(X,\Q)$
for any $X$ normal in $U$.

b) Assume conjecture $P_{-1}$. Then $H_i(X,\Q(n))=0$ for all $X$ and $n<0$,
and if $X$ has a desingularization $p:X'\to X$ which is an 
isomorphism outside of the dense open subset $U$, then 
$H_i^S(U,\Q)\cong H_i^S(X,\Q)$.
In particular, Suslin homology and higher Chow groups of
weight $0$ are birational invariant.
 
c) Under conjecture $P_0$, the groups $H_i^S(X,\Q)$ are finite dimensional
and vanish unless $0\leq i\leq d$. 

d) Conjecture $P_0$ is equivalent to the vanishing of $H_i^S(X,\Q)$ for 
all $i\not=0$ and all smooth $X$.
\end{proposition}

\proof 
The argument is the same as in Theorem \ref{ppart}. 
To prove b), we have to show
that $H^i_c(X,\Q(n))=0$ for $n>d=\dim X$ under $P_{-1}$, and for
c) we have to show that $H^i_c(X,\Q(d))$ is finite dimensional and vanishes
unless $d\leq i\leq 2d$ under $P_0$.
By induction on the dimension and the localization sequence
we can assume that $X$ is smooth and projective. In this case, the 
statement is Conjecture $P_{-1}$ and $P_0$, respectively, 
plus the fact that $H_0^S(X,\Q)\cong CH_0(X)_\Q$ is a finite dimensional 
vector space.
The final statement follows from the exact sequence \eqref{localizations}
and the vanishing of $H^i_c(X,\Q(n))=0$ for $n>d=\dim X$ under $P_{-1}$.
\proofend

\begin{proposition}
Conjecture $P_0$ holds if and only if 
the map $H^i_M(X,\Q)\to H^i_S(X,\Q)$ of \eqref{canmap}
is an isomorphism for all $X/\F_q$ and $i$.
\end{proposition}

\proof
The second statement implies the first, because if the map is an 
isomorphism, then $H^i_S(X,\Q)=0$ for $i\not=0$ and $X$ smooth and proper, 
and hence so is the dual $H_i^S(X,\Q)$.
To show that $P_0$ implies the second statement, first note that
because the map is compatible with long exact blow-up sequences,
we can by induction on the dimension assume that $X$ is smooth
of dimension $d$. In this case, motivic cohomology vanishes above
degree $0$, and the same is true for Suslin cohomology in view
of Proposition \ref{ratsituation}d). To show that for connected $X$ the map
\eqref{canmap} is an isomorphism of $\Q$ in degree zero,
we consider the commutative diagram induced by the structure map 
$$\begin{CD}
H^0_M(\F_q,\Q) @>>> H^0_S(\F_q,\Q)\\
@VVV @VVV\\
H^0_M(X,\Q) @>>> H^0_S(X,\Q) 
\end{CD}$$
This reduces the problem to the case $X=\Spec \F_q$, where it can be 
directly verified.
\proofend

\subsection{Integral coefficients}
Combining the torsion results \cite{uweshuji} with the rational results,
we obtain the following

\begin{proposition}
Conjecture $P_0$ is equivalent to the finite generation of
$H_i^S(X,\Z)$ for all $X/\F_q$.
\end{proposition}

\proof
If $X$ is smooth and proper, then according to the main theorem of 
Jannsen-Saito \cite{uweshuji}, the groups 
$H_i^S(X,\Q/\Z)=CH_0(X,i,\Q/\Z)$ are isomorphic to etale homology,
and hence finite for $i>0$ by the Weil-conjectures. 
Hence finite generation of $H_i^S(X,\Z)$ implies that $H_i^S(X,\Q)=0$.

Conversely, we can by induction on the dimension assume that $X$ is smooth
and has a smooth and proper model. Expressing Suslin homology of
smooth schemes in terms of cohomology with compact support and again
using induction, it suffices to show that $H^i_M(X,\Z(n))$ is finitely
generated for smooth and proper $X$ and $n\geq \dim X$. Using the projective
bundle formula we can assume that $n=\dim X$, and then the statement follows
because $H^i_M(X,\Z(n))\cong CH_0(X,2n-i)$ is finitely generated according
to \cite[Thm 1.1]{ichkato}.
\proofend

Recall the pairings of Lemma \ref{duall1}. We call them 
perfect if they identify one group with the
dual of the other group. In the torsion case, this implies that
the groups are finite, but in the free case this is not true:
For example, $\oplus_I \Z$ and $\prod_I \Z$ are in perfect duality.

\begin{proposition}\label{fingensushom}
Let $X$ be a separated scheme of finite type over a finite field. 
Then the following statements are equivalent:
\begin{enumerate}
\item The groups $H_i^S(X,\Z)$ are finitely generated for all $i$.
\item The groups $H^i_S(X,\Z)$ are finitely generated for all $i$.
\item The groups $H^i_S(X,\Z)$ are countable for all $i$.
\item The pairings of Lemma \ref{duall1} are perfect for all $i$.
\end{enumerate}
\end{proposition}

\proof
a) $\Rightarrow$ b) $\Rightarrow$ c) are clear, and
c) $\Rightarrow$ a) follows from \cite[Prop.3F.12]{hatcher},
which states that if $A$ is not finitely generated, then either
$\Hom(A,\Z)$ or $\Ext(A,\Z)$ is uncountable.

Going through the proof of Lemma \ref{duall1} it is easy to see
that a) implies d). Conversely, if the pairing is perfect,
then ${}_\tor H_i^S(X,\Z)$ is finite. Let $A=H^i_S(X,\Z)/\tor$
and fix a prime $l$. Then $A/l$ is a quotient of 
$H^i_S(X,\Z)/l\subseteq H^i_S(X,\Z/l)$, and which is
finite by Theorem \ref{jsfingen}.
Choose lifts $b_i\in A$ of a basis of $A/l$ and 
let $B$ be the finitely generated free abelian
subgroup of $A$ generated by the $b_i$. By construction, $A/B$ is
$l$-divisible, hence $H_i^S(X,\Z)/\tor=\Hom(A,\Z)\subseteq \Hom(B,\Z)$
is finitely generated.
\proofend

\subsection{Algebraically closure of a finite field}
Suslin homology has properties similar to a Weil-cohomology theory.
Let $X_1$ be separated and of finite type over $\F_q$, 
$X_n= X\times_{\F_q}\F_{q^n}$ and $X=X_1\times_{\F_q}\bar \F_q$.
From Corollary \ref{descentss}, we obtain a short exact sequence
$$ 0\to {\lim}^1 H^{t+1}_S(X_n,\Z)\to H^t_S(X,\Z)
\to\lim H^t_S(X_n,\Z)\to 0.$$
The outer terms can be calculated with the $6$-term $\lim$-$\lim^1$-sequence
associated to \eqref{cefff}.
The theorem of Suslin and Voevodsky implies that
$$ \lim H^i_S(X,\Z/l^r)\cong H^i_\et(X,\Z_l)$$
for $l\not=p$. For $X$ is proper and $l=p$, we get the same
result from \cite{ichduality}
$$ H^i_S(X,\Z/p^r)\cong \Hom(CH_0(X,i,Z/p^r),\Z/p^r)\cong
H^i_\et(X,\Z/p^r).$$
We show that this is true integrally:

\begin{proposition}
Let $X$ be a smooth and proper curve over the algebraic closure
of a finite field $k$ of characteristic $p$. Then the non-vanishing cohomology
groups are
$$ H^i_S(X,\Z)\cong
\begin{cases}
\Z & i=0\\
\lim_r \Hom_{GS}(\mu_{p^r}, \Pic X) \times
\prod_{l\not=p} T_l\Pic X (-1) & i=1\\
\prod_{l\not=p} \Z_l(-1) & i=2.
\end{cases} $$
The homomorphisms are maps of group schemes.
\end{proposition}

\proof
By properness and smoothness we have
$$H_i^S(X,\Z)\cong H^{2-i}_M(X,\Z(1))\cong 
\begin{cases}
\Pic X & i=0;\\
k^\times & i=1;\\
0 & i\not=0,1.
\end{cases}$$
Now
$$\Ext^1(k^\times,\Z)= 
\Hom(\colim_{p\not|m}\mu_m,\Q/\Z)\cong \prod_{l\not=p}\Z_l(-1)$$
and since $\Pic X$ is finitely generated by torsion,
\begin{multline*}
\Ext^1(\Pic X,\Z)\cong \Hom(\colim_m {}_m\Pic X,\Q/\Z)\cong \\
\lim \Hom_{GS}({}_m\Pic X,\Z/m)
\cong \lim_m \Hom_{GS}(\mu_m,{}_m\Pic X) 
\end{multline*}
by the Weil-pairing.
\proofend

\begin{proposition}
Let $X$ be smooth, projective and connected over the algebraic
closure of a finite field. Assuming conjecture $P_0$, we have
$$ H^i_S(X,\Z)\cong
\begin{cases}
\Z & i=0\\
\prod_l H^i_\et(X,\Z_l) & i\geq 1.
\end{cases} $$
In particular, the $l$-adic completion of $H^i_S(X,\Z)$
is $l$-adic cohomology $H^i_\et(X,\Z_l)$ for all $l$.
\end{proposition}

\proof
Let $d=\dim X$. By properness and smoothness we have
$$H_i^S(X,\Z)\cong H^{2d-i}_M(X,\Z(d)).$$
Under hypothesis $P_0$, the groups $H_i^S(X,\Z)$
are torsion for $i>0$, and $H_0^S(X,\Z)=CH_0(X)$ is the product of
a finitely generated group and a torsion group. Hence for $i\geq 1$
we get by \eqref{cefff} that
\begin{multline*}
H^i_S(X,\Z)\cong \Ext^1(H_{i-1}^S(X,\Z),\Z)\cong
\Hom(H_{i-1}^S(X,\Z)_\tor,\Q/\Z)\\ \cong
\Hom(H^{2d-i+1}_M(X,\Z(d))_\tor,\Q/\Z)
\cong
\Hom(H^{2d-i}_\et(X,\Q/\Z(d)),\Q/\Z) \\
\cong \Hom(\colim_m H^{2d-i}_\et(X,\Z/m(d)),\Q/\Z)
\cong
\lim_m \Hom(H^{2d-i}_\et(X,\Z/m(d)),\Z/m).
\end{multline*}
By Poincare-duality, the latter agrees with
$\lim H^i_\et(X,\Z/m) \cong \prod_l H^i_\et(X,\Z_l)$.
\proofend

\section{Arithmetic homology and cohomology}
We recall some definitions and results from \cite{ichweilI}.
Let $X$ be separated and of finite type over a finite field $\F_q$, 
$\bar X=X\times_{\F_q}\bar\F_q$ and $G$ be the Weil-group of $\F_q$.
Let $\gamma: {\cal T}_G \to {\cal T}_{\hat G}$ be the functor
from the category of $G$-modules to the category of continuous
$\hat G=\Gal(\F_q)$-modules which associated to $M$ the module
$\gamma_*=\colim_m M^{mG}$, where the index set is ordered by divisibility.
It is easy to see that the forgetful functor is a left adjoint of $\gamma_*$,
hence $\gamma_*$ is left exact and preserves limits. The derived functors
$\gamma^i_*$ vanish for $i>1$, and
$\gamma^1_*M=R^1\gamma_*M=\colim M_{mG}$, where the transition maps
are given by $M_{mG}\to M_{mnG}, x\mapsto \sum_{s\in mG/mnG}sx$.
Consequently, a complex $M^\cdot$ of $G$-modules
gives rise to an exact triangle of continuous $G_k$-modules
\begin{equation}\label{weiltriangle}
\gamma_* M^\cdot \to R\gamma_* M^\cdot \to \gamma^1_* M^\cdot[-1] .
\end{equation}
If $M=\gamma^*N$ is the restriction of a continuous $\hat G$-module, then
$\gamma_*M=N$ and $\gamma^1_*M=N\otimes \Q$. In particular, Weil-\etale
cohomology and \etale cohomology of torsion sheaves agree.
Note that the derived functors $\gamma_*$ restricted to the
category of $\hat G$-modules does not agree with the derived functors
of $\tau_*$ considered in Lemma \ref{suscohss}. Indeed, 
$R^i\tau_*M = \colim_L H^i(G_L,M)$ is the colimit of Galois cohomology
groups, whereas $R^i\gamma_*M = \colim_m H^i(mG,M)$ is the colimit
of cohomology groups of the discrete group $\Z$.

\subsection{Homology}
We define arithmetic homology with coefficients in the $G$-module $A$ to be
$$H_i^\ar(X,A) := \Tor_i^G (C_*^X(\bar k),  A).$$
A concrete representative is the double complex
$$ C_*^X(\bar k)\otimes A \stackrel{1-\varphi}{\longrightarrow}
C_*^X(\bar k)\otimes A,$$
with the left and right term in homological degrees one and 
zero, respectively, and with the Frobenius endomorphism $\varphi$ acting 
diagonally. We obtain short exact sequences
\begin{equation}\label{poip}
0 \to H_i^S(\bar X,A)_G \to H_i^\ar(X,A)\to H_{i-1}^S(\bar X,A)^G\to 0.
\end{equation}

\begin{lemma}\label{pofinite}
The groups $H_i^\ar(X,\Z/m)$ are finite. 
In particular, $H_i^\ar(X,\Z)/m$ and ${}_mH_i^\ar(X,\Z)$ are finite.
\end{lemma}

\proof
The first statement follows from the short exact sequence \eqref{poip}.
Indeed, if $m$ is prime
to the characteristic, then we apply \eqref{susvoemain} 
together with finite generation of etale cohomology, and if $m$
is a power of the characteristic, we apply Theorem \ref{ppart}
to obtain finiteness of the outer terms of \eqref{poip}.
The final statements follows from the long exact sequence
$$ \cdots \to H_i^\ar(X,\Z)\stackrel{\times m}{\longrightarrow}
 H_i^\ar(X,\Z)\to H_i^\ar(X,\Z/m)\to \cdots$$
\proofend

If $A$ is the restriction of a $\hat G$-module, then
\eqref{weiltriangle} applied to the complex of continuous
$\hat G$-modules $\C_*^X(\bar k)\otimes A$, gives a long exact sequence
$$ \cdots \to H_i^{GS}(X,A) \to  H_{i+1}^\ar(X,A) \to 
H_{i+1}^{GS}(X,A_\Q) \to H_{i-1}^{GS}(X,A) \to \cdots$$
With rational coefficients this sequence breaks up into 
isomorphisms
\begin{equation}\label{ratweilsus}
H_i^\ar(X,\Q)\cong H_i^S(X,\Q)\oplus H_{i-1}^S(X,\Q).
\end{equation}

\subsection{Cohomology}
In analogy to \eqref{suscohdef}, we define arithmetic
cohomology with coefficients in the $G$-module $A$ to be
\begin{equation}\label{weilsusdef}
H^i_\ar(X,A)= \Ext^i_G(C_*^X(\bar k),A).
\end{equation}
Note the difference to the definition in \cite{lichtenbaumweilI}, 
which does not give well-behaved (i.e. finitely generated) groups
for schemes which are not smooth and proper.
A concrete representative is the double complex
$$ \Hom(C_*^X(\bar k),A)  \stackrel{1-\varphi}{\longrightarrow}
\Hom(C_*^X(\bar k),A),$$
where the left and right hand term are in cohomological degrees
zero and one, respectively. 
There are short exact sequences
\begin{equation}\label{weilhsss}
0 \to H^{i-1}_S(\bar X,A)_G \to H^i_\ar(X,A)\to H^i_S(\bar X,A)^G\to 0.
\end{equation}

The proof of Lemma \ref{pofinite} also shows

\begin{lemma} The groups $H^i_\ar(X,\Z/m)$ are finite. 
In particular, ${}_mH^i_\ar(X,\Z)$ and $H^i_\ar(X,\Z)/m$ are finite.
\end{lemma}

\begin{lemma}
For every $G$-module $A$, we have an isomorphism
$$H^i_\ar(X,A)\cong H^i_{GS}(X,R\gamma_*\gamma^*A).$$
\end{lemma}

\proof
Since $M^G=(\gamma_*M)^{\hat G}$, Weil-Suslin cohomology is the Galois
cohomology of the
derived functor of $\gamma_*\Hom_\Ab(C_*^X(\bar k),-)$ on the 
category of $G$-modules. By Lemma \ref{suscohss}, it suffices to
show that this derived functor agrees with
the derived functor of $\tau_*\Hom_\Ab(C_*^X(\bar k),\gamma_*-)$ on the
category of $G$-modules. But given a continuous $\hat G$-modules $M$ and a
$G$-module $N$, the inclusion
$$\tau_*\Hom_\Ab(M,\gamma_*N)\subseteq  \gamma_*\Hom_\Ab(\gamma^*M,N)$$
induced by the inclusion $\gamma_*N\subseteq N$
is an isomorphism. Indeed, if $f:M\to N$ is $H$-invariant and $m\in M$
is fixed by $H'$, then $f(m)$ is fixed by $H\cap H'$, hence $f$
factors through $\gamma_*N$. 
\proofend

\begin{corollary}
If $A$ is a continuous $\hat G$-module, then there is a long exact sequence
$$ \cdots \to H^i_{GS}(X,A) \to  H^i_\ar(X,A) \to H^{i-1}_{GS}(X,A_\Q)
\to H^{i+1}_{GS}(X,A) \to \cdots.$$
\end{corollary}

\proof
This follows from the Lemma by applying the long exact
$\Ext^*_{\hat G}(C_*^X(\bar k),-)$-sequence to \eqref{weiltriangle}.
\proofend

\subsection{Finite generation and duality}

\begin{lemma}\label{duall}
There are natural pairings
$$ H^i_\ar(X,\Z)/\tor \times H_i^\ar(X,\Z)/\tor \to \Z$$
and
$$H^i_\ar(X,\Z)_\tor \times H_{i-1}^\ar(X,\Z)_\tor\to \Q/\Z.$$
\end{lemma}

\proof
From the adjunction
$ \Hom_G(M,\Z)\cong \Hom_\Ab(M_G,\Z)$
and the fact that $L(-)_G = R(-)^G[-1]$,
we obtain by deriving a quasi-isomorphism
$$ R\Hom_G(C_*^X(\bar k),\Z)\cong
R\Hom_\Ab(C_*^X(\bar k)\otimes^L_G\Z,\Z) .$$
Now we obtain the pairing as in Lemma \ref{duall1}
using the resulting spectral sequence
$$ \Ext^s_\Ab(H_t^\ar(X,\Z),\Z)\Rightarrow H^{s+t}_\ar(X,\Z).$$
\proofend

\begin{proposition}\label{fingenequi}
For a given separated scheme $X$ of finite type over $\F_q$, 
the following statements are equivalent:
\begin{enumerate}
\item The groups $H_i^\ar(X,\Z)$ are finitely generated.
\item The groups $H^i_\ar(X,\Z)$ are finitely generated.
\item The groups $H^i_\ar(X,\Z)$ are countable.
\item The pairings of Lemma \ref{duall} are perfect.
\end{enumerate}
\end{proposition}

\proof
This is proved exactly as Proposition \ref{fingensushom}, with
Theorem \ref{jsfingen} replaced by Lemma \ref{pofinite}.
\proofend

We need a Weil-version of motivic cohomology
with compact support. We define $H^i_c(X_W,\Z(n))$ to be
the $i$th cohomology of $R\Gamma(G, R\Gamma_c(\bar X,\Z(n)))$,
where the inner term is a complex defining motivic cohomology
with compact support of $\bar X$. We use this notation to distinguish
it from arithmetic homology with compact support considered in
\cite{ichduke}, which is the cohomology of 
$R\Gamma(G, R\Gamma_c(\bar X_\et,\Z(n)))$. However, if $n\geq \dim X$,
which is the case of most importance for us, both theories agree.

Similar to \eqref{localizations} we obtain for a closed subscheme
$Z$ of a smooth scheme $X$ of pure dimension $d$ with open
complement $U$ a long exact sequence
\begin{equation}\label{localizationsw}
\cdots \to H_i^\ar(U,\Z) \to H_i^\ar(X,\Z) \to
H^{2d+1-i}_c(Z_W,\Z(d))\to\cdots .
\end{equation}
The shift by $1$ in degrees occurs because arithmetic homology is defined
using homology of $G$, whereas cohomology with compact support is defined
using cohomology of $G$.

\begin{proposition}
The following statements are equivalent:
\begin{enumerate}
\item Conjecture $P_0$.
\item The groups $H_i^\ar(X,\Z)$ are finitely generated for all $X$.
\end{enumerate}
\end{proposition}

\proof
a) $\Rightarrow$ b):
By induction on the dimension of $X$ and the blow-up square, 
we can assume that $X$ is smooth of dimension $d$, where
$$H_i^\ar(X,\Z)\cong H^{2d+1-i}_c(X_W,\Z(d)).$$
By localization for $H^*_c(X_W,\Z(d))$ and induction on the dimension
we can reduce the question to $X$ smooth and projective. In 
this case $\Z(d)$ has etale hypercohomological descent
over an algebraically closed field by \cite{ichduality}, hence 
$H^j_c(X_W,\Z(d))$ agrees with the Weil-etale cohomology $H^j_W(X,\Z(d))$ 
considered in \cite{ichweilI}.
These groups are the finitely generated for $i> 2d$ 
by \cite[Thm.7.3,7.5]{ichweilI}. By conjecture $P_0$, and the
isomorphism 
$H^i_W(X,\Z(d))_\Q\cong CH_0(X,2d-i)_\Q\oplus CH_0(X,2d-i+1)_\Q$
of Thm.7.1c) loc.cit., these groups are torsion for $i<2d$, so that
the finite group $H^{i-1}(X_\et,\Q/\Z(d))$ surjects onto 
$H^i_W(X,\Z(d))$. Finally, $H^{2d}_W(X,\Z(d))$ is an extension
of the finitely generated group $CH_0(\bar X)^G$ by the finite group 
$H^{2d-1}(\bar X_\et,\Z(d))_G\cong H^{2d-2}(\bar X_\et,\Q/\Z(d))_G$.

b) $\Rightarrow$ a) Consider the special case that $X$ is smooth
and projective. Then as above, $H_i^\ar(X,\Z)\cong H^{2d+1-i}_W(X,\Z(d))$.
If this group is finitely generated, then we obtain from
the coefficient sequence that
$H^{2d+1-i}_W(X,\Z(d))\otimes \Z_l\cong \lim H^{2d+1-i}(X_\et,\Z/l^r(d))$,
and the latter group is torsion for $i>1$ by the Weil-conjectures. 
Now use \eqref{ratweilsus}.
\proofend

\begin{theorem}\label{degreezero}
For connected $X$, the map $H_0^\ar(X,\Z)\to H_0^\ar(\F_q,\Z)\cong \Z$
is an isomorphism. In particular, we have $H_0^\ar(X,\Z)\cong \Z^{\pi_0(X)}$.
\end{theorem}

\proof
The proof is similar to the proof of Theorem \ref{ratdegree0}.
Again we use induction on the dimension and the blow-up sequence to reduce 
to the situation where $X$ is irreducible and smooth. 
In this case, we can use \eqref{localizationsw}
and the following Proposition to reduce to the smooth and proper case, 
where we have
$H_0^\ar(X,\Z)=CH_0(\bar X)_G \cong \Z$.
\proofend

\begin{proposition}\label{zerohigh}
If $n>\dim X$, then $H^i_c(X_W,\Z(n))=0$
for $i>n+\dim X$. 
\end{proposition}

\proof
By induction on the dimension and the localization sequence
for motivic cohomology with compact support one sees that the
statement for $X$ and a dense open subscheme of $X$ are equivalent.
Hence we can assume that $X$ is smooth and proper of dimension $d$.
In this case, $H^i_c(X_W,\Z(n))$ is an extension of $H^i_M(\bar X,\Z(n))^G$
by $H^{i-1}_M(\bar X,\Z(n))_G$. These groups vanish for 
$i-1>d+n$ for dimension (of cycles) reasons.
For $i=d+n+1$, we have to show that $H^{d+n}_M(\bar X,\Z(n))_G$
vanishes. From the niveau spectral sequence for motivic cohomology
we obtain a surjection
$$ \bigoplus_{\bar X_{(0)}} H^{n-d}_M(k(x),\Z(n-d))
\twoheadrightarrow H^{d+n}_M(\bar X,\Z(n)).$$
The summands are isomorphic to $K_{n-d}^M(\bar\F_q)$. 
If $n>d+1$, then they vanish because higher Milnor $K$-theory of
algebraically closed fields vanishes. If $n=d+1$, then the summands 
are isomorphic to $(\bar \F_q)^\times$, whose coinvariants vanish. 
\proofend

\section{A Kato type homology}
We construct a homology theory measuring the 
difference between Suslin homology and arithmetic homology.
The cohomological theory can be defined analogously.
Kato-Suslin-homology with coefficients in the $G$-module $A$,
$H_i^{KS}(X,A)$ is defined as the $i$th homology of the complex
of coinvariants $(C_*^X(\bar k)\otimes A)_G$. If $A$ is 
trivial as a $G$-module, then since 
$(C_*^X(\bar k)\otimes A)^G \cong C_*^X(k)\otimes A$, we get the
short exact sequence of complexes
$$0\to C_*^X(k) \otimes A\to  
C_*^X(\bar k)\otimes A \stackrel{1-\varphi}{\longrightarrow}
C_*^X(\bar k)\otimes A\to (C_*^X(\bar k)\otimes A)_\varphi \to 0$$
and hence a long exact sequence
$$ \cdots \to H_i^S(X,A) \to  H_{i+1}^\ar(X,A) \to 
H_{i+1}^{KS}(X,A) \to H_{i-1}^{S}(X,A) \to \cdots .$$
By Theorem \ref{degreezero}
we have $H_0^{KS}(X,\Z)\cong H_0^\ar(X,\Z)\cong  \Z^{\pi_0(X)}$.
The following is a generalization of the integral 
version \cite{ichkato} of Kato's conjecture \cite{kato}.

\begin{conjecture}\label{katoconj} (Generalized integral Kato-conjecture)
If $X$ is smooth, then $H_i^{KS}(X,\Z)=0$ for $i>0$.
\end{conjecture}

Equivalently, the canonical map $H_i^S(X,\Z)\cong H_{i+1}^\ar(X,\Z)$
is an isomorphism for all smooth $X$ and all $i\geq 0$, i.e. there
are short exact sequences
$$0 \to H_{i+1}^S(\bar X,\Z)_G \to H_i^S(X,\Z)\to H_i^S(\bar X,\Z)^G\to 0.$$

\begin{theorem}
Conjecture \ref{katoconj} is equivalent to conjecture $P_0$.
\end{theorem}

\proof
If Conjecture \ref{katoconj} holds, then 
$$ H_i^S(X,\Q)\cong H_{i+1}^\ar(X,\Q)\cong H_{i+1}^S(X,\Q)\oplus
H_i^S(X,\Q)$$
implies the vanishing of $H_i^S(X,\Q)$ for $i>0$.

Conversely, we first claim that for smooth and proper $Z$, 
the canonical map $H^i_c(Z,\Z(n))\to H^i_c(Z_W,\Z(n))$ is an isomorphism for 
all $i$ if $n<\dim Z$, and for $i\leq 2n$ if $n=\dim Z$.  
Indeed, if $n\geq \dim Z$ then the cohomology of $\Z(n)$ agrees with the 
etale hypercohomology of $\Z(n)$, see \cite{ichduality},
hence satisfies Galois descent. 
But according to (the proof of) Proposition \ref{ratsituation}b), 
these groups are torsion groups, so that Galois descent $R\Gamma_{G_k}$
agrees with $R\Gamma_G$.

Using localization for cohomology with compact support and 
induction on the dimension, we get next that 
$H^i_c(Z,\Z(n))\cong H^i_c(Z_W,\Z(n))$ for all $i$ and all $Z$ 
with $n<\dim Z$. Now choose a smooth
and proper compactification $C$ of $X$.  Comparing the 
exact sequences \eqref{localizations} and \eqref{localizationsw}, 
we see with the $5$-Lemma that the 
isomorphism $H_i^S(C,\Z)\cong H^{2d-i}_c(C,\Z(d))
\to H_{i+1}^\ar(C,\Z)\cong H^{2d-i}_c(C_W,\Z(d))$ for $C$ implies the 
same isomorphism for $X$ and $i\geq 0$.
\proofend

\section{Tamely ramified class field theory}
We propose the following conjecture relating Weil-Suslin homology 
to class field theory:

\begin{conjecture}\label{CFT} (Tame reciprocity conjecture)
For any $X$ separated and of finite
type over a finite field, there is a canonical injection to 
the tame abelianized fundamental group with dense image 
$$ H_1^\ar(X,\Z) \to \pi_1^{t}(X)^{ab}.$$
\end{conjecture}

Note that the group $H_1^\ar(X,\Z)$ is conjecturally finitely generated.
At this point, we do not have an explicit construction (associating
elements in the Galois groups to algebraic cycles) of the map.
One might even hope that 
$H_1^\ar(X,\Z)^\circ := \ker (H_1^\ar(X,\Z)\to \Z^{\pi_0(X)}$
is finitely generated and isomorphic to the abelianized 
geometric part of the tame fundamental group defined in SGA 3X\S 6.

Under Conjecture \ref{katoconj}, $H_0^S(X,\Z)\cong H_1^\ar(X,\Z)$ 
for smooth $X$, 
and Conjecture \ref{CFT} is a theorem of Schmidt-Spiess \cite{ss}.

\begin{proposition}
a) We have $H_1^\ar(X,\Z)^{\wedge l}\cong \pi_1^t(X)^{ab}(l)$.
In particular, the prime to $p$-part of Conjecture \ref{CFT} holds
if $H_1^\ar(X,\Z)$ is finitely generated. 

b) The analog statement holds for the $p$-part if 
$X$ has a compactification $T$ which has a
desingularization which is an isomorphism outside of $X$. 
\end{proposition}

\proof
a) By Theorem \ref{degreezero}, $H_0^\ar(X,\Z)$ contains
no divisible subgroup. Hence if $l\not=p$, we have by 
Theorems \ref{svduality} and \ref{ffduality}
\begin{multline*}
H_1^\ar(X,\Z)^{\wedge l} \cong \lim H_1^\ar(X,\Z/l^r)\cong 
\lim H_0^{GS}(X,\Z/l^r) \\
\cong \lim H^1_\et(X,\Z/l^r)^*
\cong \pi_1^t(X)^{ab}(l).
\end{multline*}

b) Under the above hypothesis, we can use the duality result
of \cite{ichduality} for the proper scheme $T$ to get 
with Proposition \ref{biratt}
\begin{multline*} 
H_1^\ar(X,\Z)\otimes \Z_p\cong  \lim H_0^{GS}(X,\Z/p^r)\cong
\lim H_0^{GS}(T,\Z/p^r)\\
\cong \lim H^1_\et(T,\Z/p^r)^*
\cong \pi_1(T)^{ab}(p)\cong \pi_1^t(X)^{ab}(p).
\end{multline*}
\proofend

\end{document}